\documentclass[11pt]{article}

\title{A Bijection Between Two Different Classes of Partitions Enumerated by $p_\nu(n)$}
\author{A.S. Andersen}
\usepackage[utf8]{inputenc}
\usepackage{biblatex}
\usepackage{youngtab}
\usepackage{amsmath}  
\usepackage{amsfonts} 
\usepackage{amssymb}
\newtheorem{theorem}{Theorem}

\begin{document}
\maketitle
\centerline{\bf Abstract}
In this paper, we give a purely bijective proof that two different partition classes that are both combinatorial interpretations of the partition function $p_\nu(n)$, a partition function related to the third order mock theta function $\nu(q)$, are equinumerous. In doing so, we give a partial solution to a combinatorial problem proposed in a paper by Andrews.
\section{Introduction and Notation}
Consider the third order mock theta function $\nu(q)$, which was first defined by Watson \cite{4} and may be defined as follows:
\begin{equation}
\nu(q) := \sum\limits_{n=0}^{\infty}\frac{q^{n(n+1)}}{(-q;q^2)_{n+1}} ,
\end{equation}
where the q-Pochhammer symbol $(a;q)_n$ is defined as usual
\begin{equation}
(a;q)_n := \prod\limits_{k=0}^{n-1}(1-aq^k) .
\end{equation}
The partition function $p_\nu(n)$ may be defined as the partition function for which $\nu(-q)$ is the generating function, and a number of combinatorial interpretations have been given for this partition function. Among these is the number of self-conjugate odd Ferrers graphs of 2$n$+1 and the number of self-conjugate partitions of 4$n$+1 into odd parts \cite{2}, \cite{3}. Odd Ferrers graphs, introduced by Andrews in \cite{1}, may be defined as Ferrers graphs in which a 2 is placed in every box, except the surrounding border, where 1s are placed in each box. For example, the following odd Ferrers graph represents the partition 7 + 7 + 3 + 1:
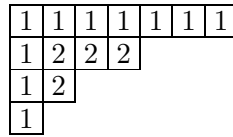
\[\young(1111111,1222,12,1)\]. Let $\mathcal{O}_{2n+1}$ be the set of self-conjugate odd Ferrers graphs for 2$n$+1, let $\mathcal{S}_{4n+1}$ be the set of self-conjugate partitions of 4$n$+1 into odd parts, let $\mathcal{O} = \cup_{n > 0}\mathcal{O}_{2n+1}$ and let $\mathcal{S} = \cup_{n > 0}\mathcal{S}_{4n+1}$. The following theorem has previously been proven through non-bijective means \cite{2} :  
\begin{theorem}
$|\mathcal{O}_{2n+1}| = |\mathcal{S}_{4n+1}|$ for all $n$.
\end{theorem}
We will give a purely bijective proof of this theorem by describing a bijection $\phi$ such that $\phi(\lambda) = \mu$, where $\lambda$ and $\mu$ are both partitions, $\lambda \in \mathcal{O}_{2n+1}$, and $\mu \in \mathcal{S}_{4n+1}$, and use the case where $\lambda$ = 3 + 5 + 3, representable as the following odd Ferrers graph:
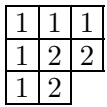
\[\young(111,122,12)\] 
as an example (Note that in this example case, $\lambda \in \mathcal{O}_{11}$, and that $\mu \in \mathcal{S}_{21}$). In doing so,  we give a partial solution to the combinatorial challenge proposed by Andrews \cite{2} asking for bijections between the various classes of partitions enumerated by $p_\nu(n)$.
\section{A Bijection Between $\mathcal{O}_{2n+1}$ and $\mathcal{S}_{4n+1}$}
 Consider the fact that the Ferrers diagrams of self-conjugate partitions may be thought of as being made up of "hooks" of other self-conjugate partitions in which every part other than the greatest part is equal to 1. For example, the Ferrers digram of the self-conjugate partition 4 + 4 + 2 + 2
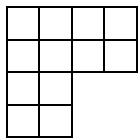
\[\young(~~~~,~~~~,~~,~~)\]
can be thought of as consisting of the following "hooks":
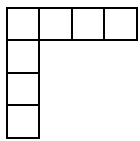
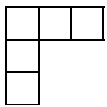
\[\young(~~~~,~,~,~)\] and \[\young(~~~,~,~)\].
Let $h_i$ denote the $i$th "hook" in a self-conjugate partition $\pi$, where $i > 0$. Note that, where $|\pi|$ may denote the sum of the parts of $\pi$, where $|h_i|$ may denote the sum of the parts in each hook in the Ferrers digram of $\pi$, and where $n$ may denote the number of hooks in $\pi$, that $\sum\limits_{i=1}^{n}|h_i| = |\pi|$. Additionally, for $\lambda \in \mathcal{O}$, let $t = \sum\limits_{i=2}^{n}|h_i|$, or the sum of the 2s in the odd Ferrers diagram. We will distinguish between the hooks in $\lambda$ and the hooks in $\mu$ by using $h_i$ to denote the $i$th hook in the former and $\eta_i$ to denote the $i$th hook in the latter. The map $\phi(\lambda) = \mu$ may be described as follows:
\newline
\linebreak
\linebreak
Step 1: Create $\eta_1$ by creating a hook with the largest part equal to $|h_1|$. Note that $|\eta_1| = 2|h_1| - 1$. For the example case for $\lambda$ given above, $\eta_1$ would be the following:
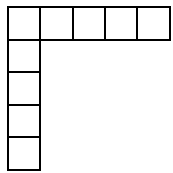
\[\young(~~~~~,~,~,~,~)\]
Step 2: For each $h_i$ where $i > 1$, create $\eta_{2i-2}$ such that $|\eta_{2i-2}| = |h_i|+1$, and $\eta_{2i-1}$ such that $|\eta_{2i-1}| = |h_i|-1$. For example, in the example case of $\lambda$ given above, $|h_2|$ = 6, so we create $\eta_2$ and $\eta_3$ such that $|\eta_2| = 7$ and $|\eta_3| = 5$, and since the number of hooks in $\lambda$ is equal to 2, the creation of these hooks completes the bijection resulting in the following partition:
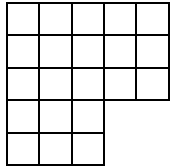
\[\young(~~~~~,~~~~~,~~~~~,~~~,~~~)\]
or 5 + 5 + 5 + 3 + 3. The map described evidently always results in a self-conjugate partition. The map described also always results in a partition of 4$n$+1, because in creating $\eta_1$ we create a partition of size $2h_1 - 1$, and in adding every $\eta_i$ such that $i >1$, we add $2t$ to this partition, thus making a partition of size $2(h_1 + t) - 1$. We know that $h_1 + t = |\lambda| = 2n + 1$, so substituting $2n + 1$ for $h_1 + t$ in the previous expression reveals that the sum of the parts in the newly created partition is always equal to $4n + 1$. Additionally, we know that the newly created partition is always a partition into odd parts because it always creates a partition in which the greatest part of $\eta_1$ is odd, the number of hooks is odd, and in which the greatest part of each hook alternates in parity, where the greatest part of $\eta_{2i-2}$ is always one greater than the greatest part of $\eta_{2i-1}$. The inverse map is obvious, so $\phi$ is a bijection, and thus $|\mathcal{O}_{2n+1}| = |\mathcal{S}_{4n+1}|$ for all $n$.
\section{Further Remarks}
Recall the natural bijection that exists between the class of self-conjugate partitions of $n$ and the class of partitions of $n$ into distinct odd parts that maps a self-conjugate partition onto a partition into distinct odd parts by making the sum of the parts in each of the hooks in the self conjugate partition into a part in the newly created partition. Where $\mathcal{D}_{2n+1}$ may denote the set of partitions of 2$n$ + 1 into distinct parts in which there is 1 odd part which is greater than half the greatest even part and every other part is even and is of the form 4$k$ + 2 where $k \in \mathbb{N}$, and where $\mathcal{DO}_{4n+1}$ may denote the set of partitions of 4$n$ + 1 into an odd number of distinct odd parts such that, when ordered from largest to smallest, the parts alternate between being of the form 4$k$ + 1 and being of the form 4$k$ + 3 where again $k \in \mathbb{N}$, an analogous bijection exists between $\mathcal{O}_{2n+1}$ and $\mathcal{D}_{2n+1}$ and between $\mathcal{S}_{4n+1}$ and $\mathcal{DO}_{4n+1}$. Thus, the bijection given above induces one between $\mathcal{D}_{2n+1}$ and $\mathcal{DO}_{4n+1}$.
\section*{Acknowledgements}
The author would like to thank George Andrews and Shane Chern for their helpful comments and suggestions.

\end{document}